\documentclass[12pt]{amsart}
\usepackage{amsthm}
\usepackage{amsfonts}
\usepackage{amsmath}
\usepackage{amssymb}
\usepackage{graphicx}
\usepackage{float} 
\newtheorem{theorem}{Theorem}[section]
\newtheorem{corollary}{Corollary}[theorem]
\newtheorem{lemma}[theorem]{Lemma}
\newtheorem{proposition}{Proposition}
\newtheorem{prop}{Proposition}
\newtheorem{example}{Example}

\newtheorem{definition}{Definition}
\newtheorem{remark}{Remark}

\begin{document}


\oddsidemargin 16.5mm
\evensidemargin 16.5mm

\thispagestyle{plain}

\vspace{5cc}
\begin{center}

{\large\bf  Watson-Crick strong bi-catenation on words
\rule{0mm}{6mm}\renewcommand{\thefootnote}{}
\footnotetext{\scriptsize ${}^{\ast}$Corresponding author. Kalpana Mahalingam}
\footnotetext{\scriptsize 2020 Mathematics Subject Classification. FILL SUBJECT MSCs HERE.

\rule{2.4mm}{0mm}Keywords and Phrases. Binary Operation, Bi-Catenation, Watson-Crick powers}}

\vspace{1cc}
{\large\it Kalpana Mahalingam \\
Department Of Mathematics, \\
Indian Institute of Technology, \\
Chennai-600036. India.\\
Email~:~kmahalingam@iitm.ac.in}

\vspace{1cc}
\parbox{24cc}{{\small

In this paper we  define and investigate the binary word operation of strong-$\phi$-bi-catenation (denoted by $\leftrightarrows_\phi$) where $\phi$ is either a morphic or an antimorphic involution. In particular, we concentrate on the mapping $\phi=\theta_{DNA}$, which  models the Watson-Crick complementarity of DNA single strands.  We show that such an operation is commutative and not associative and when iteratively applied to a word $u$, this operation generates words over $\{u, \theta(u)\}$. We then extend this operation to languages and show that the families of regular, context-free and context-sensitive languages
are closed under the operation of strong-$\phi$-bi-catenation. We also define the notion of $\leftrightarrows_\theta$-conjugacy and study conditions on words $u$ and $v$ where $u$ is a $\leftrightarrows_\theta$-conjugate of $v$. We then extend this relation to language equations  and provide solutions under some special cases.}}
\end{center}


\vspace{1.5cc}
\section{Introduction}\label{secintro}
Combinatorics on words focuses on the study of words and formal languages(\cite{Lothaire97,ssyu2005}). A word is basically formed from alphabets by simply juxtaposing the alphabets. Such an operation is called as concatenation, which is indeed a basic binary operation on words. 
Some of the well known basic word operations defined and studied in literature are quotient, shuffle(\cite{Eilen,Lshuf}), bi-catenation(\cite{Bcat}), $k$-catenation(\cite{kcat}), insertion(\cite{Lthesis}) and deletion(\cite{Lthesis}) to name a few. These operations were naturally extended to languages and authors in general studied closure properties of the families in the Chomsky hierarchy under the above operations among others. 

In \cite{kinvcat}, the authors used the $k$-catenation operation defined in \cite{kcat} to define $k$-involution codes. The $k$-involution codes formally denote DNA strands (possibly used in DNA based computations) avoiding certain non-specific hybridizations that pose potential problems for the results of the biocomputation.  A DNA strand is basically a word over the alphabet $\{A,G,C,T\}$ and its Watson-Crick complement is mathematically formalized by an antimorphic involution denoted by $\theta_{DNA}$ which is an antimorphism ($\theta_{DNA}(uv) =\theta_{DNA}(v)\theta_{DNA}(u)$) and an involution ($\theta^2_{DNA}(u) = u$) that maps $A\mapsto T$, $C\mapsto G$ and vice-versa. Concatenation of DNA strands is a process to combine various DNA strands linearly to form new DNA strands. One such recombination is obtained by repeatedly concatenating a DNA strand $u$ and its Watson-Crick complement $\theta_{DNA}(u)$ in random order. Such a strand is called a $\theta_{DNA}$-power of $u$. The authors in \cite{man14,cai2023} extended the notion of catenation to $\phi$-catenation and strong $\phi$-catenation respectively that generates all possible $\phi$-powers of a given word where $\phi$ is either a morphic or an antimorphic involution.

 Observe that, the operation strong-$\phi$-catenation, when applied iteratively to a word $u$, results in all possible $\phi$-powers of $u$ (i.e.) words that belong to the set $\{u,\phi(u)\}^+$. However, when this operation is applied between two distinct words, say $u$ and $v$, the resulting set does not provide all possible combinations of words of the set $\{u,v, \phi(u),\phi(v)\}$ as the catenation is one-sided. To fill this gap, we introduce the notion of strong-$\phi$-bi-catenation of words. 
 
 In this paper, we combine the notion of strong-$\phi$-catenation and bi-catenation to obtain a new binary operation which we call strong-$\phi$-bi-catenation of words. We
define and investigate some basic properties of strong-$\phi$-bi-catenation in Section \ref{defpro}. We also mention its connection to the
previously defined notion of strong-$\phi$-catenation. In Section \ref{Lext}, we naturally extend the operation to languages and show that the families of regular, context-free and context-sensitive languages are closed under this operation.

Section \ref{Cop} briefly explore closure properties of languages closed under strong-$\phi$-bi-catenation. Section \ref{concom} investigates conjugacy with respect to $\leftrightarrows_\phi$ and Section \ref{equations} studies some language equations with respect to the strong $\phi$-bi-catenation operation.  We end the paper with few concluding remarks.

\section{Preliminaries}
\label{secprelim}

Let $ \Sigma$ be a finite alphabet. We denote by $\Sigma ^*$ the set of all words over $\Sigma$ including the empty word $\lambda$. By $\Sigma^+$, we denote the set of all non-empty words over $\Sigma$. The length of a word $w \in \Sigma^*$ is the number of letter occurrences in $w$, denoted by $|w|$; i.e. if $ w = a_1 a_2 ... a_n,~ a_i \in \Sigma $ then $|w|=n$. $|w|_a$ denotes the number of occurrences of $a$ in $w$.
The reverse of the word $w=a_1a_2\cdots a_{n-1}a_n$ denoted by $w^R$ is the word $a_na_{n-1}\cdots a_2a_1$ where $a_i \in \Sigma$, $1\le i \le n$.
A word $w$ is called primitive it is not the non-trivial power of another word; i.e. if $w = u^i$ then $w=u$ and $i=1$. The \textit{primitive root} of a word $w$ is the shortest $u$ such that $w=u^i$ for some $i$, denoted by $\rho(w)=u$. We denote by $\mathcal{Q}$, the set of all primitive words.

We first recall some results from \cite{lyndon,cai2023}.
\begin{lemma}\cite{lyndon} \label{lemconj}
Let $u,\;v,\;w \in \Sigma^{+}$ be such that,
 $uv=vw$, then for $k \geq 0$, $x \in \Sigma^{+}$ and $y \in \Sigma^{*}$, $u=xy$, $v=(xy)^{k}x$, $w=yx$.
 \end{lemma}

\begin{lemma}\label{lpow}\cite{lyndon}
If $xy=yx$ then $x$ and $y$ are powers of a common word; i.e. $x=u^i$ and $ y=u^j$ for some $u\in \Sigma^+$.
\end{lemma}

\begin{lemma}\label{reflem7}\cite{cai2023}
For $x,y \in \Sigma^+$, if $yxx =xxy$, then $x =\alpha^m$ and $y =\alpha^n$  for some $m,n\ge 1$ and $\alpha \in \Sigma^+$.
\end{lemma}

A mapping $\phi:\Sigma^{*} \rightarrow \Sigma^{*}$ is called  a \textit{morphism} on $\Sigma^*$ if for all words $u,v \in \Sigma^{*}$ we have that $\phi(uv)=\phi(u) \phi(v)$, an \textit{antimorphism} on $\Sigma^*$ if $\phi(uv)=\phi(v) \phi(u)$ and an \textit{involution} if $\phi(\phi(a)) = a$ for all $a\in \Sigma$.

A mapping $\phi: \Sigma^* \rightarrow \Sigma^*$ is called a {\it morphic involution on $\Sigma^*$}  (respectively, an {\it antimorphic involution on $\Sigma^*$}) if it is an involution  on $\Sigma$ extended to a morphism (respectively, to an antimorphism) on $\Sigma^*$.
For convenience, in the remainder of this paper  we use the convention that the letter  $\phi$   denotes an involution that is either morphic or antimorphic (such a mapping will be termed {\it (anti)morphic involution}),  that  the letter  $\theta$   denotes an antimorphic involution, and that the letter  $\mu$  denotes a morphic involution.
For $L\subseteq \Sigma^*$ and an involution $\phi$, we define,
\[\phi(L) = \{\phi(w)~:~w\in L\}\]
\[L^R = \{w^R~:~ w \in L\}.\]
A word $u$ is a conjugate of $v$ if for some $w$, $uw=wv$. Two words $u$ and $v$ are said to commute if $uv = vu$. The concept of conjugacy and commutativity was extended to the notion of an involution map $\theta$ in \cite{Kari08}. Recall that $u$ is said to be a $\theta$-conjugate of $w$ if $uv=\theta(v)w$ for some $v \in \Sigma^+$, and  $u$ is said to $\theta$-commute with $v$ if $uv = \theta(v) u$. We recall the following result  from \cite{Kari08} characterizing $\theta$-conjugacy and $\theta$-commutativity for an antimorphic involution $\theta$ (if $\theta = \theta_{DNA}$,  these are called Watson-Crick conjugacy, respectively Watson-Crick commutativity). For an antimorphic involution $\theta$, a word $u$ is called a $\theta$-palindrome if $u=\theta(u)$. The set of all $\theta$-palindromes is denoted by $P_\theta$.

\begin{proposition}\label{gg1}\cite{Kari08} For $u,v,w \in \Sigma^+$ and
 $\theta$ an antimorphic involution,
 \begin{enumerate}
     \item If $uv = \theta(v)w$, then  either there exists $x\in \Sigma^+$ and $y\in \Sigma^*$ such that $u =xy $ and $w=y\theta(x)$, or $u =\theta(w)$.
 \item If $uv = \theta(v)u$, then $u = x(yx)^i$, $v = yx$, for some $i \geq 0$ and   $\theta$-palindromes $x \in \Sigma^*, y\in \Sigma^+$.
\end{enumerate}
\end{proposition}

We recall the following from \cite{Kari2010}.
\begin{proposition}\cite{Kari2010}\label{palpro1}
Let $x, y \in  \Sigma^+$ and
 $\theta$ an antimorphic involution, such that $xy = \theta(y)\theta(x)$ and $yx = \theta(x)\theta(y)$. Then, one of the following holds:
\begin{enumerate}
    \item $x=\alpha^i$, $y =\alpha^k$ for some $\alpha \in P_{\theta}$ 
    \item $x=[\theta(s)s]^i\theta(s)$, $y = [s\theta(s)]^ks$ for some $s\in \Sigma^+$, $i,k\ge 0$.
\end{enumerate}
\end{proposition}





We recall the following from \cite{cai2023}.
\begin{definition} For a given $u \in \Sigma^*$, and an (anti)morphic involution $\phi$, the set $\{u,\phi(u)\}$ is denoted by $u_{\phi}$,  and is called a $\phi$-complementary pair, or $\phi$-pair for short.  The length of a $\phi$-pair  $u_{\phi}$ is defined as  $|u_{\phi}| = |u| = |\phi(u)|$.
 \end{definition}

It was also remarked in \cite{cai2023} that, for $u\in \Sigma^+$ and $\phi$, an (anti)morphic involution,
 $|u_{\phi}|_a = |u|_a + |\phi(u)|_a$,
$|\phi(u)|_a = |u|_{\phi(a)}$ and $|\phi(u)|_{\phi(a)} = |u|_a$. 
 For $L\subseteq\Sigma^*$, we denote   $L_\phi = L\cup \phi(L)$. A word is called $\phi$-power of a word $u$ if it is of the form $u_1 u_2 ... u_n$ where $u_1 = u$ and $u_i \in u_\phi$ for $2\leq i \leq n $. 

\section{Strong $\phi$-bi-catenation}\label{defpro}

In this section, we define and study a new binary operation called the strong $\phi$-bi-catenation. The basic string operation \textit{catenation} is a binary operation that maps $(u,v)$ to $uv$. The catenation operation has several generalizations. The first one is the notion of \textit{Bi-catenation} (\cite{Bcat}), which is a binary operation which maps $(u,v)$ to $\{uv, vu\}$. Motivated by the Watson-crick complemantarity of DNA strands, the authors in \cite{man14},  defined the concept of \textit{$\phi$-catenation} which incorporates an (anti)-morphic involution mapping $\phi$. The $\phi$-catenation maps $(u,v)$ to $\{uv,u\phi(v)\}$. This concept was further generalized in \cite{cai2023} to define the \textit{strong $\phi$-catenation}, which generates all possible $\phi$ powers of a given word $u$, (i.e.) all words in the set $\{u,\phi(u)\}^+$. In this section, we introduce the notion of \textit{strong $\phi$-bi-catenation} operation which is indeed a generalization of bi-catenation defined in \cite{Bcat} as well as the strong $\phi$-catenation operation(\cite{cai2023}).

Binary operation $\circ$ on $\Sigma^*$ is a map $\circ : \Sigma^* \times \Sigma^* \to 2^{\Sigma^*}$. For a given binary operation $\circ$, the {\it  $i$-th $\circ$-power of a word} is defined by :
    $$u^{\circ(0)} = \{\lambda\}, ~u^{\circ(1)} = u \circ \lambda, ~u^{\circ(i)} = u^{\circ(i-1)}\circ u, ~i\ge 2$$

Note that, depending on the operation $\circ$, the $i$-th power of a word can be a singleton word, or a set of words.

A binary operation called $\phi$-catenation denoted by $\odot$, was defined in \cite{man14} which generates some $\phi$ powers of a word $u$ under consideration, when $\odot$ is applied iteratively. However, this concept was extended to the notion of strong-$\phi$-catenation denoted by $\otimes$,
  that generates all the non-trivial  $\phi$-powers of $u$, that is, the union of the sets $\{u, \theta(u)\}^n$, $n\geq 2$.  We begin the section by recalling the formal definition of strong $\phi$-catenation.

\begin{definition}\cite{cai2023}\label{defstro}
Given an (anti)morphic involution $\phi$ on $\Sigma^*$ and  two words $u,v\in \Sigma^*$, we define the strong-$\phi$-catenation operation  of $u$ and $v$ with respect to $\phi$  as
$$ u\otimes v = \{uv, u\phi(v), \phi(u)v, \phi(u)\phi(v)\}.$$
\end{definition}
We recall the following from \cite{cai2023}.
\begin{proposition}\label{omr1}\cite{cai2023}
For an antimorphic involution $\theta$ and $u,v\in \Sigma^+$, $u\otimes v = v\otimes u$ iff {\it(i)} $u=v$, or {\it (ii)} $u =\theta(v)$, or {\it (iii)} $u$ and $v$ are powers of a common $\theta$-palindrome.
\end{proposition}

We now formally define the notion of strong $\phi$-bi-catenation operation.

\begin{definition}
We define strong $\phi$-bi-catenation  $({\leftrightarrows_{\phi}})$ as 
$$u \leftrightarrows_{\phi} v = (u\otimes v)\cup (v\otimes u) = u_\phi v_\phi  \cup v_\phi u_\phi$$
Writing explicitly all the terms of $u\leftrightarrows_{\phi} v$ we get,  \[
     u \leftrightarrows_{\phi} v = \{uv, u\phi(v), \phi(u)v, \phi(u)\phi(v), vu , v\phi(u), \phi(v)u, \phi(v)\phi(u)\} 
\]

\begin{example}\label{exn}
Consider  the case of $\theta =\theta_{DNA}$, the Watson-Crick complementary function that maps $A\leftrightarrow T$ and $C\leftrightarrow G$  and the words $ u= ATC$, $v= GCTA$. Then,
\begin{align*}
    u \leftrightarrows_{\theta} v &= 
\{ATCGCTA,~ ATCTAGC, ~GATGCTA,~GATTAGC \} \\
& \cup \{GCTAATC, ~GCTAGAT, ~TAGCATC, ~TAGCGAT \}
\end{align*}
which  is the set of all bi-catenations that involve words $u$ and $v$ and their images under $\theta_{DNA}$.
\end{example}


\end{definition}

We have the following remark which follows directly from definition.
\begin{remark}\label{rem1}
     Let $\phi$ be an (anti)morphic involution on $\Sigma^*$ and $u,v \in \Sigma^*$. Then, for $u_1\in u_\phi$ and $v_1\in v_\phi$,
     $$u\leftrightarrows_{\phi}v = u_1\leftrightarrows_{\phi}v_1 = v_1\leftrightarrows_{\phi}u_1$$
\end{remark}

We first observe the following which is straightforward from the definition.

\begin{lemma}
    Let $\phi$ be an (anti)morphic involution on $\Sigma^*$ and $u,v \in \Sigma^*$. Then,  $x \in u\leftrightarrows_{\phi}v$ iff $\phi(x) \in u\leftrightarrows_{\phi}v$.
\end{lemma}



 
 
 A bw-operation $\circ$ is called length-increasing if for any $u, v \in \Sigma^+$ and $w \in u \circ v$,
$|w| > max\{|u|, |v|\}$.
A bw-operation $\circ$ is called propagating if for any $u, v \in \Sigma^*$, $a \in \Sigma$ and $w \in u \circ v$,
$|w|_a = |u|_a + |v|_a$. In \cite{man14}, these notions were generalized to incorporate an (anti)morphic involution $\phi$, as follows. A bw-operation $\circ$ is called $\phi$-propagating if
for any $u, v \in \Sigma^*$, $a \in \Sigma$ and $w \in u \circ v$, $|w|_{a,\phi(a)} = |u|_{a,\phi(a)}+|v|_{a,\phi(a)}$. 
It was shown in \cite{man14} that the operation $\phi$-catenation is  not propagating but is $\phi$-propagating. The concept of $\phi$-catenation was extended to strong $\phi$-catenation in \cite{cai2023}. It was shown in \cite{cai2023} that the operation strong $\phi$-catenation is also not propagating but is $\phi$-propagating. 

A bw-operation $\circ$ is called left-inclusive if for any three words  $u, v, w \in \Sigma^*$ we have
$$(u \circ v) \circ w \supseteq  u \circ (v \circ w)$$
and is called right-inclusive if
$$(u \circ v) \circ w \subseteq u \circ (v \circ w).$$

A bw-operation $\circ$ is associative if for any three words $u,v,w \in \Sigma^*$ we have 
$$(u\circ v)\circ w = u\circ (v\circ w)$$
Similar to the properties  of the operation $\phi$-catenation and strong $\phi$-catenation investigated in \cite{man14,cai2023}, one can easily observe that the strong-$\phi$-bi-catenation operation is  length increasing, not propagating and $\phi$-propagating.  In \cite{man14},  it was shown that for a  morphic  involution the  $\phi$-catenation operation is  trivially  associative, whereas for  an antimorphic involution the  $\phi$-catenation  operation is not associative. In contrast, it was shown in \cite{cai2023}, that the strong-$\phi$-catenation operation is right inclusive, left inclusive, as well as associative, when $\phi$ is a morphic as well as an antimorphic involution. We also observe that the operation strong $\phi$-bi-catenation operation is commutative and not associative.

\begin{lemma}
Let $\phi$ be an (anti)morphic  involution. The strong $\phi$-bi-catenation operation is length increasing, not propagating,  $\phi$-propagating, commutative, not associative, and neither right nor left inclusive.
\end{lemma}
\begin{proof}
We show that the binary operation $\leftrightarrows_{\phi}$ is length increasing, $\phi$-propagating and commutative.
\begin{enumerate}
    \item Let $u,v,w \in \Sigma^+$ such that $w \in u\leftrightarrows_{\phi} v$. Then, $|w| = |u| +|v|$ and hence $|w| > max\{|u|, |v|\}$. Thus, the operation $\leftrightarrows_{\phi}$ is length increasing.
    \item Consider the words $u,v$ from Example \ref{exn}. Note that, for $w = GATGCTA \in u\leftrightarrows_{\phi} v$, $|w|_G = 2 \neq |u|_G +|v|_G = 0+1 =1$. Hence, the operation $\leftrightarrows_{\phi}$ is not propagating.

    \item Let $w, u,v\in\Sigma^+$ be such that $w \in u\leftrightarrows_{\phi} v$. Then, $$w \in \{uv, u\phi(v), \phi(u)v, \phi(u)\phi(v), vu , v\phi(u), \phi(v)u, \phi(v)\phi(u)\}$$
 Suppose, $w = \phi(v)u$ then,
 \begin{equation*}
    \begin{split}
    |w|_{a,\phi(a)} & = |w|_a + |w|_{\phi(a)} \\
    &=|\phi(v)|_{a,\phi(a)} + |u|_{a,\phi(a)} \\
    &= |\phi(v)|_{\phi(a)} +|\phi(v)|_a  + |u|_{a, \phi(a)}  \\
    &= |v|_{a} +|v|_{\phi(a)}  + |u|_{a, \phi(a)}  \\
    &= |u|_{a,\phi(a)} + |v|_{a,\phi(a)}
    \end{split}
\end{equation*}
The other cases are similar and we omit them.  Hence, the operation $\leftrightarrows_{\phi}$ is $\phi$-propagating.
\item One can easily observe from the definition that for $u,v\in \Sigma^*$,
$$u\leftrightarrows_\phi v = u_\phi v_\phi\cup v_\phi u_\phi = v\leftrightarrows_\phi u$$
Hence, $\leftrightarrows_\phi$ is commutative.
    \item 
       Note that, for $u =AG$, $v=CA$ and $w=AC$ and $\theta =\theta_{DNA}$, we have $CACTAC \in  v_{\phi}u_{\phi}w_{\phi} \subseteq (u\leftrightarrows_{\phi} v)\leftrightarrows_{\phi} w $ but not in $ u\leftrightarrows_{\phi} (v\leftrightarrows_{\phi} w) $
     Thus, the operation $\leftrightarrows_{\phi}$ is not associative.
     \item It is evident from the example given in Item 5 that the operation $\leftrightarrows_{\phi}$ is neither right nor left inclusive. \qed
 \end{enumerate}   
\end{proof}

We now give a sufficient condition on  words $u$ and $w$ such that  
$(u\leftrightarrows_{\phi} v)\leftrightarrows_{\phi} w = u\leftrightarrows_{\phi} (v\leftrightarrows_{\phi} w)$.
\begin{lemma}
    Given an (anti)morphic involution $\phi$ and $u,v,w\in \Sigma^+$ such that $u_\phi w_\phi = w_\phi u_\phi$ then, 
    $$(u\leftrightarrows_{\phi} v)\leftrightarrows_{\phi} w = u\leftrightarrows_{\phi} (v\leftrightarrows_{\phi} w).$$
\end{lemma}
\begin{proof}
    Let $u,v,w\in \Sigma^+$. Then,
    \begin{align*}
        (u\leftrightarrows_{\phi} v)\leftrightarrows_{\phi} w &= \{u_{\phi}v_{\phi}~\cup~ v_{\phi}u_{\phi}\} \leftrightarrows_{\phi} w \\
&= u_{\phi}v_{\phi}w_{\phi} \cup v_{\phi}u_{\phi}w_{\phi} \cup w_{\phi}u_{\phi}v_{\phi} \cup w_{\phi}v_{\phi}u_{\phi}
     \end{align*}    
     and,
      \begin{align*}
        u\leftrightarrows_{\phi} (v\leftrightarrows_{\phi} w) &= u\leftrightarrows_{\phi}\{v_{\phi}w_{\phi}~\cup~ w_{\phi}v_{\phi}\}  \\
&= u_{\phi}v_{\phi}w_{\phi} \cup u_{\phi}w_{\phi}v_{\phi} \cup v_{\phi}w_{\phi}u_{\phi} \cup w_{\phi}v_{\phi}u_{\phi}
     \end{align*} 
     Thus, if $u_\phi w_\phi = w_\phi u_\phi$, then $(u\leftrightarrows_{\phi} v)\leftrightarrows_{\phi} w = u\leftrightarrows_{\phi} (v\leftrightarrows_{\phi} w)$.
\end{proof}

 \subsection{Extension to Languages}\label{Lext}
 In this section we extend the $\leftrightarrows_{\phi}$ operation to languages. We use the notation $L_{\phi}$ to denote the set $L\cup \phi(L)$.
Given $L_1,L_2\subseteq \Sigma^*$ define, 
$$L_1 \leftrightarrows_{\phi} L_2 =
 \bigcup_{u\in L_1, v\in L_2}
(u\leftrightarrows_{\phi} v)$$ and 
$L_1 \leftrightarrows_{\phi} \emptyset = \emptyset \leftrightarrows_{\phi} L_2 =\emptyset$ and $L_1 \leftrightarrows_{\phi}^{0}  L_2 =L_1\cup \phi(L_1)\cup L_2\cup \phi(L_2)$. The iterated strong -bi-$\phi$-catenation operation $\leftrightarrows_{\phi}^ i$ for $i\ge 1$ and languages $L_1$ and $L_2$ is defined as
$L_1 \leftrightarrows_{\phi}^i L_2 = (L_1 \leftrightarrows_{\phi}^{i-1} L_2) \leftrightarrows_{\phi} L_2$. 
The $i$-th $\leftrightarrows_{\phi}$-power of a non-empty
language $L$  is defined as
$$L^{\leftrightarrows_{\phi}(0)} = \{\lambda\}, ~ L^{\leftrightarrows_{\phi}(1)} = L_{\phi},~ L^{\leftrightarrows_{\phi}(i)}= (L \leftrightarrows_{\phi}^{i-1} L),~ i\ge 1$$ 
The $+$-closure of a non-empty language $L$ with respect to a bw-operation $\leftrightarrows_{\phi}$, denoted
by $L^{\leftrightarrows_{\phi}(+)}$ is defined as $$L^{\leftrightarrows_{\phi}(+)} =\bigcup_{k\ge 1} L^{\leftrightarrows_{\phi}(k)}$$ 
\\
 We say that $L$ is $\leftrightarrows_{\phi}$-\textit{closed} if for any $u$ and $v$ in $L$, $u \leftrightarrows_{\phi} v$ is a subset of $L$.  We say that a binary operation $\leftrightarrows_{\phi}$ is $plus$-$closed$ if for any non-empty language $L \subset \Sigma^*$,   $L^{\leftrightarrows_{\phi}(+)}$ is also $\leftrightarrows_{\phi}$-closed.
\\

We first observe that, 
$u \leftrightarrows_{\phi} u = u\otimes u$ and hence, $u^{\leftrightarrows_\phi(n)} = u^{\otimes(n)}$ for all $n\ge 0$.
Thus, for $u=ATC$ and $\theta=\theta_{DNA}$ we have,
$$u^{\leftrightarrows_\phi(n)}= u^{\otimes(n)} = \{u_1u_2\cdots u_n~:~ u_i = ATC ~\text{or}~ u_i = GAT, ~1\le i\le n\}$$


We observe the following.
\begin{lemma}\label{lemch}
    For a language $U,V\subset \Sigma^* $, $$U \leftrightarrows_{\phi}V= U_{\phi}V_{\phi}~ \cup ~V_{\phi}U_{\phi} $$
\end{lemma}
\begin{proof}
For $U,V\subseteq \Sigma^*$, we have,
    \begin{equation*}
    \begin{split}
    U\leftrightarrows_{\phi}V &= \bigcup_{u_1\in U, u_2 \in V}  u_1 \leftrightarrows_{\phi} u_2 \\
    &= \bigcup_{u_1\in U, u_2 \in V}  ( (u_1)_\phi(u_2)_\phi \cup (u_2)_\phi(u_1)_\phi) )\\
    & = \bigcup_{u_1\in U, u_2 \in V}  (u_1\otimes u_2) \cup (u_2\otimes u_1) \\
    & = U\otimes V \cup V\otimes U\\
    & = U_{\phi}V_{\phi}~ \cup ~V_{\phi}U_{\phi} \\
\end{split}
\end{equation*}
\end{proof}

We now have the following observation which characterizes the form of words
in $L^{\leftrightarrows_{\phi}(n)}$ when the strong-$\phi$-bi-catenation operation is applied iteratively.
\begin{prop}\label{pro1212}
For a language $ L \subset \Sigma^* $, $L^{\leftrightarrows_{\phi}(n)}$ is the collection of all words of the form $u_1 u_2 ... u_n$ where $u_i \in L_\phi$ and  $n \geq 2$.
\end{prop}
\begin{proof}
We use induction on $n$. For $n=2$,
\begin{equation*}
    \begin{split}
    L^{\leftrightarrows_{\phi}(2)} = L \leftrightarrows_{\phi} L &= \bigcup_{u_1, u_2 \in L}  u_1 \leftrightarrows_{\phi} u_2 \\
    &= \bigcup_{u_1, u_2 \in L} ( (u_1)_\phi(u_2)_\phi \cup (u_2)_\phi(u_1)_\phi) )\\
    & = \bigcup_{u_1, u_2 \in L}  (u_1)_\phi(u_2)_\phi\\
    & = \{ u_1 u_2~:~u_1 , u_2 \in L_\phi \} 
    \end{split}
\end{equation*}
Now assume that $ L^{\leftrightarrows_{\phi}(n)} = \{ u_1 u_2 ... u_n ~:~ u_i \in L_\phi\}  $. For $n+1$, 
\begin{equation*}
    \begin{split}
    L^{\leftrightarrows_{\phi}(n+1)} = L^{\leftrightarrows_{\phi}(n)} \leftrightarrows_{\phi} L &= \bigcup_{u\in L^{\leftrightarrows_{\phi}(n)}, u' \in L}  u \leftrightarrows_{\phi} u' \\
    &= \bigcup_{u\in L^{\leftrightarrows_{\phi}(n)}, u' \in L} ( (u)_\phi(u')_\phi \cup (u')_\theta(u)_\theta) )\\
    & = \bigcup_{u\in L^{\leftrightarrows_{\phi}(n)}, u' \in L}  (u)_\phi(u')_\phi\\
    & = \{ u u' ~:~ u \in L^{\leftrightarrows_{\phi}(n)} , u' \in L_\phi \} \\
    & = \{ u_1 u_2 ... u_n u_{n+1}~:~ u_i \in L_\phi \}
    \end{split}
\end{equation*}
Hence the result.
\end{proof}

\begin{proposition}\label{pro1213}
Let $ L \subset \Sigma^* $. For any morphic or antimorphic involution, $$L^{\leftrightarrows_{\phi}(n)} \leftrightarrows_{\phi} L^{\leftrightarrows_{\phi}(m)} = L^{\leftrightarrows_{\phi}(n+m)}$$
\end{proposition} 
\begin{proof}
Using the above result (Proposition \ref{pro1212}), we have 
\begin{align*}
  L^{\leftrightarrows_{\phi}(n+1)} &= L^{\leftrightarrows_{\phi}(n)} \leftrightarrows_{\phi} L^{\leftrightarrows_{\phi}(1)}\\
    &= \{ u_1 u_2 ... u_n u_{n+1}~:~ u_i \in L_\phi \}
\end{align*} 
Repeating the $\leftrightarrows_{\phi}$ operation $m$ times and using above result (Proposition \ref{pro1212}) we have, \[L^{\leftrightarrows_{\phi}(n+m)} = \{ u_1 u_2 ... u_{n+m}~:~ u_i \in L_\phi \} = L^{\leftrightarrows_{\phi}(n)} \leftrightarrows_{\phi} L^{\leftrightarrows_{\phi}(m)}\]
\end{proof}

\begin{corollary}\label{corp} The operation $\leftrightarrows_{\phi}$ is plus-closed; i.e., for any $ u , v \in L^{\leftrightarrows_{\phi}(+)} $, we have $u \leftrightarrows_{\phi} v \in L^{\leftrightarrows_{\phi}(+)}$.
\end{corollary}

\begin{proof}
Let $ u , v \in L^{\leftrightarrows_{\phi}(+)}$. Then, there exist $n$ and $m$ such that $u \in L^{\leftrightarrows_{\phi}(n)}$ and $v \in L^{\leftrightarrows_{\phi}(m)}$. By Proposition \ref{pro1213}, we have $u \leftrightarrows_{\phi} v \in L^{\leftrightarrows_{\phi}(n+m)} $. Thus, $u \leftrightarrows_{\phi} v \in L^{\leftrightarrows_{\phi}(+)} $.
\end{proof}

One can also easily observe that for a regular (context-free, context-sensitive) language $L$, $\phi(L)$ is also regular (context-free, context-sensitive respectively). Thus, from Lemma \ref{lemch}, we conclude the following.
\begin{theorem}
    The families of regular, context-free and context-sensitive languages
are closed under the operation of strong bi-$\phi$-catenation.
\end{theorem}

\subsection{$\leftrightarrows_\theta$-closed Languages}\label{Cop}
A language $L$ is closed under the mapping $\phi$ if  $x \in L$ implies $\phi(x) \in L$ i.e., $L= L_{\phi}$ and is closed under catenation if $u,v\in L$, imply $uv\in L$.
A language $L$ is $\leftrightarrows_\phi$-closed if $u, v \in L$ imply
$u \leftrightarrows_\phi v \subseteq L$. It was shown in Corollary \ref{corp}  that the operation $\leftrightarrows_\phi$ is plus-closed. 

\begin{lemma}\label{lemclo1}
    If $L$ is closed under $\phi$ and catenation then $L$ is closed under $\leftrightarrows_\phi$.
\end{lemma}
\begin{proof}
    If $L$ is closed under $\phi$ then $L = L_{\phi}$ and if $L$ is closed under catenation then, $L^2 =L$. From Lemma \ref{lemch} we observe that, $L\leftrightarrows_\phi L = L_{\phi}L_{\phi} = L^2 = L$. Hence, $L$ is closed under $\leftrightarrows_\phi$.
    \end{proof}
    The converse of Lemma \ref{lemclo1} is not true in general. For example, consider the alphabet $\{a,b\}$ and an antimorphic involution $\theta$ such that $\theta(a)=a$ and $\theta(b)=b$. Let $L = \{ ab\} \cup \{x~:~ x\in \{a,b\}^+,~|x|\ge 3\}$. Note that, $L$ is closed under catenation and $L$ is closed under $\leftrightarrows_\theta$ but $L$ is not closed under $\theta$ as $\theta(ab) = ba \notin L$. \\
    
    We now give an example of a language $L$ such that $L$ is closed under $\leftrightarrows_\phi$.
    \begin{example}\label{exc-1}
    Consider the alphabet $\{a,b\}$ and $\phi$ be an (anti)morphic involution that maps $a$ to $b$ and vice-versa. Let $L = \{w~:~ |w|_a=|w|_b \} \subseteq \Sigma^+$. Note that for any $x \in L$, $\phi(x) \in L$ and for $x,y \in L$, $xy \in L$. Hence by Lemma \ref{lemclo1}, $L$ is closed under $\leftrightarrows_\phi$.
\end{example}

 \begin{lemma}\label{lemssc}
        Let $L$ be such that $L$ is  $\leftrightarrows_\phi$ closed. Then, $L_1L_2L_3\cdots L_n \subseteq L$ for $L_i \in L_{\phi}$ for $1\le i\le n$ and $n\ge 2$.
\end{lemma}
 \begin{proof}
     We first observe from Lemma \ref{lemch} that, \[L\leftrightarrows_\phi L = L_\phi L_\phi = L^2 \cup L\phi(L) \cup \phi(L)L \cup \phi(L)\phi(L).\] Since, $L$ is closed under $\leftrightarrows_\phi$, we have $L\leftrightarrows_\phi L \subseteq L$ which implies that $L_1L_2 \subseteq L$ for $L_1, L_2\in L_\phi$. 
     One can easily prove by induction that,  $L^n\leftrightarrows_\phi L^n = L_1L_2L_3\cdots L_n$ for $L_i \in L_{\phi}$ for $1\le i\le n$ and $n\ge 2$. Since $L_1L_2 \subseteq L$ for $L_1, L_2\in L_\phi$ we have that $L_1L_2L_3\cdots L_n\subseteq L$ for $L_i \in L_{\phi}$ for $1\le i\le n$ and $n\ge 2$ and hence the result.
 \end{proof}       

    \begin{lemma}
        Let $L$ be such that $L$ is  $\leftrightarrows_\phi$ closed. Then, the following are true.
        \begin{enumerate}
            \item $L$ is closed under catenation.
            \item $L^R$ is closed under $\leftrightarrows_\phi$.
            \item $\phi(L)$ is closed under $\leftrightarrows_\phi$.
           
            \item For all $A,B\in L_{\phi}^n$, $A \leftrightarrows_\phi B \subseteq L$.
            \end{enumerate}
    \end{lemma} 
    \begin{proof}Given that $L$ is closed under $\leftrightarrows_\phi$ (i.e.) for all $u,v \in L$, we have $u\leftrightarrows_\phi v \subseteq L$. Then let, 
    \begin{align*}
        A&= u\leftrightarrows_\phi v  \\&= \{uv, \phi(u)v, u\phi(v), \phi(u)\phi(v), vu,  \phi(v)u, v\phi(u), \phi(v)\phi(u) \}\\
        &\subseteq L
    \end{align*} 
    \begin{enumerate}
        \item  Note that,  $u\leftrightarrows_\phi v \subseteq L$ implies that $uv \in L$ for all $u,v \in L$. Hence, $L$ is closed under catenation.
        \item For $u,v \in L$ we have $u^R,v^R\in L^R$. Then, 
        $A= u\leftrightarrows_\phi v \subseteq L$ and 
        $A^R = \{u^Rv^R, \phi(u^R)v^R, u^R\phi(v^R), \phi(u^R)\phi(v^R), v^Ru^R, \phi(v^R)u^R, \\ v^R\phi(u^R), \phi(v^R)\phi(u^R) \}= u^R\leftrightarrows_\phi v^R\subseteq L^R$. Hence, $L^R$ is closed under $\leftrightarrows_\phi$.
        \item It is easy to observe that,  $u\leftrightarrows_\phi v = \phi(u)\leftrightarrows_\phi \phi(v)$ and $A= u\leftrightarrows_\phi v \subseteq L$ implies $\phi(A) \subseteq \phi(L)$. But, $A= u\leftrightarrows_\phi v = \phi(u)\leftrightarrows_\phi \phi(v) = \phi(A) \subseteq \phi(L)$. Thus, $\phi(L)$ is  closed under $\leftrightarrows_\phi$.
       
        \item Since $L$ is closed under $\leftrightarrows_\phi$, we have by Lemma \ref{lemssc}, $A \leftrightarrows_\phi B = L^n \leftrightarrows_\phi L^n = L_1L_2L_3\cdots L_n\subseteq L$  for all $n\ge 2$ and $A,B\in L_{\phi}^n$, $L_i \in \{L, \phi(L)\}$. Hence, the result.
    \end{enumerate}  
\end{proof} 
We now have the following example.
\begin{example}\label{excl-1}
    Consider the alphabet $\{a,b,c\}$ and $\phi$ an (anti)morphic involution  that maps $a$ to $b$ and vice-versa and $\phi(c) = c$. Let $L_1 = \{w~:~ |w|_a+|w|_b = |w|_c \}$ and $L_2 = \{ w ~:~ |w|_a = |w|_b =|w|_c\}$.  Note that for any $x \in L_1$, $\phi(x) \in L_1$ and for $x,y \in L_1$, $xy \in L_1$. Hence, $L_1$ is closed under $\leftrightarrows_\phi$. Similarly one can verify that $L_2$ is closed under $\leftrightarrows_\phi$. 
\end{example}
It is clear from the above example that in general for a given $\leftrightarrows_\phi$-closed language $L_1$, $L_1^c$ is not $\leftrightarrows_\phi$-closed. 

\begin{lemma}
Let $L_1,L_2\subseteq \Sigma^+$ be such that $L_1$ and $L_2$ are closed under $\leftrightarrows_\phi$. Then the following are true. 
    \begin{enumerate}
     \item $L_1^c$ is not closed under $\leftrightarrows_\phi$.
        \item $L_1\cap L_2$ is  closed under $\leftrightarrows_\phi$.
        \item $L_1\cup L_2$ is not closed under $\leftrightarrows_\phi$.
        
     \end{enumerate}
\end{lemma}
\begin{proof}
\begin{enumerate}
 \item Consider the language $L_1=  \{w~:~ |w|_a=|w|_b \}$ discussed in Example \ref{exc-1}. Then, $L_1^c =  \{w~:~ |w|_a\neq |w|_b \}$ and for $u= aba, v=bab\in L_1^c$, we have $uv= ababab \in u\leftrightarrows_\phi v$ but $uv\notin L_1^c$. Thus, for a given $L$ which is closed under $\leftrightarrows_\phi$, $L_1^c$ is not necessarily closed under $\leftrightarrows_\phi$.
    \item Given that $L_1$ and $L_2$ are closed under $\leftrightarrows_\phi$. Let $u,v \in L_1\cap L_2$. Then, $u\leftrightarrows_\phi v \in L_1\cap L_2$. Thus, $L_1\cap L_2$ is closed under $\leftrightarrows_\phi$.
    \item Consider $L_1$ and $L_2$  from Example \ref{excl-1}. Note that, $L_1 =\phi(L_1)$, $L_2 =\phi(L_2)$ and $ abc, bcca \in L_1\cup L_2$. But, $abcbcca \in abc \leftrightarrows_\phi bcca \nsubseteq L_1 \cup L_2$.  Hence, $L_1\cup L_2$ is not $\leftrightarrows_\phi$-closed.
\end{enumerate}
    \end{proof}

We now define the $\leftrightarrows_\phi$-Iterative closure of a language $L$ denoted by $cl_{\leftrightarrows_\phi}(L)$
 \begin{definition}\label{defcl}
     For a given language $L\subseteq \Sigma^+$, we define the  $\leftrightarrows_\phi$-Iterative closure of a language $L$ denoted by $cl_{\leftrightarrows_\phi}(L) = \bigcup_{i\ge 0} L_i$ where $L_0 = L_{\phi}$,  $$L_i = \{u\leftrightarrows_\phi v ~:~ u,v\in \bigcup_{k=0}^{i-1}L_{k}\}.$$
 \end{definition}
 We have the following observation which is clear from Definition \ref{defcl}.
 \begin{lemma}
   For $L\subseteq \Sigma^*$, 
   $$cl_{\leftrightarrows_\phi}(L) = \{x_1x_2\cdots x_n~:~ n\ge 1, x_i \in L_{\phi}\}= L_{\phi}^{(+)}$$
 \end{lemma}
     Note that for each $i\ge 0$, $L_i$ defined above is $\phi$-closed. Also, observe that the iterative closure of a language $L$, denoted by $cl_{\leftrightarrows_\phi}(L)$ is $\leftrightarrows_\phi$-closed. 
     
\begin{example}
     Consider the alphabet $\{a,b\}$ and $\theta$ an antimorphic involution  such that $\theta(a)=b$ and vice-versa. Let $L = \{ab\}$. Note that, $L_{\theta} = L = L_0$. Then, $L_1 = \{abab\}$, $L_2 = \{(ab)^2, (ab)^3, (ab)^4\}$ and $L_n = \{ (ab)^i~:~ 2 \le i\le 2n\}$. Hence, $cl_{\leftrightarrows_\theta}(L)  = \{(ab)^i~:~ i\ge 1\}$
 \end{example}
 \begin{theorem}
     The families of regular, context-free and context sensitive languages are closed  under the iterative $\leftrightarrows_\phi$-closure operation.
 \end{theorem}


\section{Conjugacy of words with respect  to $\leftrightarrows_\phi$}
\label{concom}
The conjugate of a word is one of the basic concept in combinatorics of words. A word $u$ is called a conjugate of $v$ if both $u$ and $v$ satisfy the word equation $u\boldsymbol{\cdot} w=w\boldsymbol{\cdot} v$ for some word $w\in \Sigma^*$ where $\boldsymbol{\cdot}$ represents the basic catenation operation. This catenation operation can be replaced by any binary operation $\circ$ to define a $\circ$- conjugate of a given word (i.e.) $u$ is a $\circ$-conjugate of $v$, if there exists a $w\in \Sigma^*$ such that $u\circ w = w\circ v$. Depending on the operation $\circ$, $u\circ w$ may be a singleton or a set. The authors in \cite{cai2023}, studied properties of $u$ and $v$  when $u$ is a $\otimes$-conjugate of $v$. 
\\
In this section, we discuss conditions  on words $u,w\in \Sigma^+$, such that $u$ is a $\leftrightarrows_\phi$-conjugate of $w$, i.e., $u \leftrightarrows_\phi v = v\leftrightarrows_\phi w$ for some $v \in \Sigma^+$. The special case when $u=w$ always holds true by definition, as the operation $\leftrightarrows_\phi$ is commutative. Thus, we can say that $u$ $\leftrightarrows_\phi$-commutes with $v$ for all $u,v \in \Sigma^*$. We prove a necessary and sufficient condition for $\leftrightarrows_\phi$-conjugacy (Theorem \ref{mr1}). Since the Watson-Crick complementarity function  $\theta_{DNA}$ is an antimorphic involution, in the remainder of this paper we only investigate  antimorphic involution mappings $\phi =\theta$.

\begin{proposition}\label{pcj1}
 Let $u, v,w\in \Sigma^+$ be such that $uv =vw$ and $u\leftrightarrows_\theta v = v\leftrightarrows_\theta w$. Then, one of the following hold true.
 \begin{enumerate}
          \item $u = s^m = w$ and $v= s^n$ for some $s\in \Sigma^+$.
     \item  $u=p^m =w$ and $v =p^n$, for $p \in P_{\theta}$.
     \item $u=xy$, $v=(xy)^ix$ and $w=yx$ for $x,y\in P_{\theta}$ and $i\ge 0$.
 \end{enumerate}
\end{proposition}
\begin{proof}
By definition, for $u,v, w\in \Sigma^+$, 
$$u\leftrightarrows_\theta v = \{uv, vu,  u\theta(v), v\theta(u), \theta(u)v, \theta(v)u, \theta(u)\theta(v), \theta(v)\theta(u)\}$$ and similarly,
$$v\leftrightarrows_\theta w = \{vw, wv, v\theta(w), w\theta(v), \theta(v)w, \theta(w)v, \theta(v)\theta(w), \theta(w)\theta(v)\}$$
Given that $uv=vw$ and $u\leftrightarrows_\theta v = v\leftrightarrows_\theta w$. Then, by Lemma \ref{lemconj}, we have $u= xy$, $v=(xy)^ix$ and $w=yx$. 
We now have the following cases.
\begin{enumerate}
    \item If $u\theta(v) = v\theta(w)$ then, $u\theta(v) = (xy)(\theta(x)\theta(y))^i\theta(x) = (xy)^ix\theta(x)\theta(y)$. If $i\neq 0$ then, $x,y \in P_{\theta}$ and $xy=yx$ and hence, $u,v$ and $w$ are powers of a common $\theta$-palindrome. If $i=0$ then, $xy\theta(x) =x\theta(x)\theta(y)$ and by Proposition \ref{gg1}, $y=st$ and $\theta(x) =(st)^js$ where $s,t \in P_{\theta}$ and hence, $x\in P_{\theta}$. Thus, $$u\leftrightarrows_\theta v = \{xyx, \theta(y)xx, xxy,x\theta(y)x \}$$
    and $$v \leftrightarrows_\theta w = \{xyx, xx\theta(y), yxx, x\theta(y)x\}$$
    Since, $u\leftrightarrows_\theta v =v\leftrightarrows_\theta w$, we have either $\theta(y)xx = xx\theta(y)$ or $\theta(y)xx = yxx$. If $\theta(y)xx = xx\theta(y)$ then, by Lemma \ref{reflem7}, $x = p^{m_1}$, $y= p^{m_2}$ for $p\in P_{\theta}$. Thus, $u = p^m=w$, $v=p^n$ for $p\in P_{\theta}$. If $\theta(y)xx = yxx$  then, $y\in P_\theta$ which implies that $u=xy$, $v=x$ and $w=yx$ for $x,y\in P_{\theta}$.
    \item The case when  $ u\theta(v) =\theta(v)w$ 
    is similar to case (1) and we omit it.
\item If $u\theta(v) = \theta(v)\theta(w)$ then, $u\theta(v) = xy(\theta(x)\theta(y))^i\theta(x) =$ $(\theta(x)\theta(y))^i\theta(x)
\theta(x)\theta(y)$ $=\theta(v)\theta(w)$. If $i=0$ then, $x\in P_{\theta}$ and the case is similar to the previous one.
If $i\neq 0$ then, $x,y \in P_{\theta}$ and $yx = xy$ and hence, $y=p^{j_1}$, $x = p^{j_2}$. Thus, $u=w =p^m$ and $v= p^n$ for $p\in P_{\theta}$.
\item If $u\theta(v) = wv$ then, $v=\theta(v)$ and $u=w$ which implies that $u=xy = w=yx$ which implies that $x$ and $y$ are powers of a common word.  If $i=0$ then, $v = x \in P_{\theta}$ and 
$$u\leftrightarrows_\theta v = \{xyx, xxy, x\theta(y)x, \theta(y)xx\}$$
$$v\leftrightarrows_\theta w = \{xyx, yxx, xx\theta(y), x\theta(y)x\}$$
and the case is similar to the previous one. If $i\neq 0$ then, $v=\theta(v)$ implies that $(xy)^ix = (\theta(x)\theta(y)^i\theta(x)$ which implies that $x,y \in P_{\theta}$. Thus, in both cases we get, $x$ and $y$ to be powers of a common $\theta$-palindromic word. 

\item If $u\theta(v) = w\theta(v)$ then, $u=w$ which implies that $u=xy =w=yx$. Hence, $x$ and $y$ are powers of a common word. Thus, $u = s^m = w$ and $v= s^n$ for some $s\in \Sigma^+$.
Then,
$$u\leftrightarrows_\theta v = \{s^k, s^m\theta(s^n), s^n\theta(s^m), \theta(s^m)s^n,  \theta(s^k)\} =v\leftrightarrows_\theta w $$
\item If $u\theta(v) = \theta(w)v$ then, $u=\theta(w) = xy = \theta(x)\theta(y)$ and $v=\theta(v)$.
Thus, $u =xy$, $v=(xy)^ix$ and $w = yx$ where $x,y \in P_{\theta}.$
\item The case when $u\theta(v) = \theta(w)\theta(v)$ is similar to the previous case and we omit it. 
\end{enumerate}    
\end{proof}

A similar proof works for the next result and hence, we omit it.
\begin{proposition}\label{pcj2}
 Let $u, v,w\in \Sigma^+$ be such that $uv =v\theta(w)$ and $u\leftrightarrows_\theta  v = v\leftrightarrows_\theta  w$. Then, one of the following hold true.
 \begin{enumerate}
     \item $u =\theta(w)=(pq)^{j+1}p$, $v= (pq)^jp$ for some $p,q\in \Sigma^*$ and $j\ge 0$.
     \item $u =w=(pq)^{j+1}p$, $v= (pq)^jp$ for some $p,q\in P_{\theta}$ and $j\ge 0$.
     \item $u=\alpha^m =w$ and $v =\alpha^n$, for $\alpha\in P_{\theta}$.
     \item $u=w = xy$ and $v=(xy)^ix$, for $x,y \in P_{\theta}$ and $i\ge 0$.
 \end{enumerate}
\end{proposition}
We now have the following result which is used in Proposition \ref{pcj3}.
\begin{lemma}\label{lus1}
    Let $x,y \in \Sigma^+$ be such that $xx\theta(y) =\theta(y)\theta(x)x$ for an antimorphic involution $\theta$.
    Then, $x$ and $y$ are powers of a common $\theta$-palindromic word.
\end{lemma}
\begin{proof}
    Given that $xx\theta(y) = \theta(y)\theta(x)x$. Then by Lemma \ref{lemconj}, we have $xx=pq$, $\theta(y) = (pq)^ip$ and $\theta(x)x = qp$ for some $p,q\in \Sigma^+$. If $|x| \le |p|$ then, $x = p_1=p_2q$ where $p=p_1p_2$. Then, $\theta(x)x = \theta(p_2q)p_1 = qp$ which implies that $q\in P_{\theta}$ and $p_1p_2 = \theta(p_2)p_1$.  By Lemma \ref{lemconj}, we have, $p_1 =\alpha(\beta\alpha)^j$ and $p_2 =\beta\alpha$ for some $\alpha, \beta \in P_{\theta}$, $j\ge 0$. Thus, $x = \alpha(\beta\alpha)^j= p_2q = \beta\alpha q$ which implies that $j\neq 0$ and $\alpha \beta = \beta\alpha$. Hence, by Lemma \ref{lpow}, $\alpha$ and $\beta$ are powers of a common word. Therefore, $x$ and $y$ are powers of a common $\theta$-palindromic word.
    The case when $|x|\ge |p|$ is similar and we omit it.
\end{proof}
\begin{proposition}\label{pcj3}
 Let $u, v,w\in \Sigma^+$ be such that $uv =\theta(v)w$ and $u\leftrightarrows_\theta  v = v\leftrightarrows_\theta  w$. Then, one of the following hold true.
 \begin{enumerate}
 \item $u=\theta(w)$ and $v=\gamma w$ for some $\gamma \in P_{\theta}$.
     \item  $u = (xy)^{j+1}p=w$ and $v= yx$ for $x,y\in P_{\theta}$ and $j\ge 0$.
     \item $u=xy=\theta(w)$, $v= \theta(x)$ for $y\in P_{\theta}$.
     \item $u=xy=\theta(w)$, $v= x$ for $x, y\in P_{\theta}$.
     \item $u=\alpha^m =w$ and $v =\alpha^n$, for $\alpha\in P_{\theta}$.
     \item $u = \theta(t)s^{k}=\theta(w)$ and $v=s^{n}t $  where $s=t\theta(t)$.
    
 \end{enumerate}
\end{proposition}
\begin{proof} 
By definition, for $u,v, w\in \Sigma^+$, 
$$u\leftrightarrows_\theta v = \{uv, vu,  u\theta(v), v\theta(u), \theta(u)v, \theta(v)u, \theta(u)\theta(v), \theta(v)\theta(u)\}$$ and similarly,
$$v\leftrightarrows_\theta w = \{vw, wv, v\theta(w), w\theta(v), \theta(v)w, \theta(w)v, \theta(v)\theta(w), \theta(w)\theta(v)\}$$
    Given that  $uv =\theta(v)w$ and $u\leftrightarrows_\theta  v = v\leftrightarrows_\theta  w$. Then by Proposition \ref{gg1}, we have either $u=\theta(w)$ and $v=\gamma w$ for some $\gamma \in P_{\theta}$ or $u=xy$, $v=\theta(x)$ , $w=y\theta(x)$ for some $x,y\in \Sigma^*$.  
    If $u=xy$, $v=\theta(x)$ , $w=y\theta(x)$ for some $x,y\in \Sigma^*$ then,
\begin{align*}
u\leftrightarrows_\theta v &= \{xy\theta(x), \theta(x)xy, xyx, \theta(x)\theta(y)\theta(x), \theta(y)\theta(x)\theta(x), xxy, \\
&~~~~~~\theta(y)\theta(x)x, x\theta(y)\theta(x)\}
\end{align*}
and similarly,
 \begin{align*}
 v\leftrightarrows_\theta w &= \{xy\theta(x), \theta(x)y\theta(x), \theta(x)x\theta(y),  xx\theta(y), y\theta(x)\theta(x),\\
 &~~~~~~y\theta(x)x, x\theta(y)\theta(x), x\theta(y)x\}
 \end{align*}
 We have the following subcases.
 \begin{enumerate}
     \item If $\theta(x)xy = \theta(x)y\theta(x)$ then, $xy = y\theta(x)$ and by Lemma \ref{lemconj}, $x=pq$, $y=(pq)^jp$ for some $p,q\in P_{\theta}$ which implies that,
     $u = (pq)^{j+1}p=w$ and $v= qp$.
     \item If $\theta(x)xy = \theta(x)x\theta(y)$ then, $y\in P_{\theta}$. Thus, $u=xy$, $v= \theta(x)$, $w=y\theta(x)$ for $y\in P_{\theta}$.
     \item If $\theta(x)xy = xx\theta(y)$ then, $u=xy$, $v=x$ and $w =yx$ for $x,y\in P_{\theta}$.
     \item If $\theta(x)xy = y\theta(x)\theta(x)$ then by Lemma \ref{lus1}, $x$ and $y$ are powers of a common $\theta$-palindromic word and hence, $u,v$ and $w$ are powers of a common $\theta$-palindromic word.
     \item If $\theta(x)xy = y\theta(x)x$ then by Lemma \ref{lpow}, $y$ and $\theta(x)x$ are powers of a common word say $s$. Then, $y = s^m$ and $\theta(x)x = s^n$. If $\theta(x) = s^{n_1}$ then, $s \in P_\theta$ and $u,v$ and $w$ are powers of a $\theta$-palindromic word $s$. If $\theta(x) = s^{n_1}s_1$ then, $x= s_2s^{n_1}$ for $s=s_1s_2$ and $2n_1 +1 =n$ which implies that $s_1 = \theta(s_2)=t$ and $s =t\theta(t)$. Thus, $u = \theta(t)s^{n_1+m}$, $v=s^{n_1}t $ and $w = s^{m+n_1}t$.
    \item If $\theta(x)xy = x\theta(y)\theta(x)$ then, $x \in P_{\theta}$ and $xy = \theta(y)x$ and by Lemma \ref{lemconj}, we have, $\theta(y) = pq$, $x=(pq)^jp$ where $p,q\in P_\theta$. Thus, $u = (pq)^{j+1}p$, $v =(pq)^jp $ and $w= qp(pq)^jp$. Therefore, $u\leftrightarrows_\theta v = v\leftrightarrows_\phi w$ implies that $pq =qp$ and by Lemma \ref{lpow}, $p$, $q$ and hence, $u,v$ and $w$ are powers of a common $\theta$-palindromic word.
     \item 
     The case when $\theta(x)xy = x\theta(y)x$ is similar to the previous and we obtain $u,v$ and $w$ to be powers of a common $\theta$-palindromic word.
    
 \end{enumerate}
\end{proof}
The proof of the following is similar to that of Proposition \ref{pcj3}.

\begin{proposition}\label{pcj4}
 Let $u, v,w\in \Sigma^+$ be such that $uv =\theta(v)\theta(w)$ and $u\leftrightarrows_\theta  v = v\leftrightarrows_\theta  w$. Then, one of the following hold true.
 \begin{enumerate}
 \item $u=w$ and $v=\gamma \theta(w)$ for some $\gamma \in P_{\theta}$.
     \item  $u = (xy)^{j+1}p=w$ and $v= yx$ for $x,y\in P_{\theta}$ and $j\ge 0$.
     \item $u=xy$, $v= \theta(x)$, $w=y\theta(x)$ for $y\in P_{\theta}$.
     \item $u=xy$, $v= x$, $w=yx$ for $x, y\in P_{\theta}$.
     \item $u=\alpha^m =w$ and $v =\alpha^n$, for $\alpha\in P_{\theta}$.
     \item $u = \theta(t)s^{k}$, $v=s^{n}t $ and $w = s^{k}t$ where $s=t\theta(t)$.
     
 \end{enumerate}
 \end{proposition}

 Based on the above results (Propositions \ref{pcj1},~\ref{pcj2},~\ref{pcj3} and \ref{pcj4}), we give a necessary and sufficient condition on words $u$, $v$ and $w$ such that $u\leftrightarrows_\theta v = v\leftrightarrows_\theta w$.

\begin{theorem}\label{mr1}
 Let $u, v,w\in \Sigma^+$. Then, $u\leftrightarrows_\theta v = v\leftrightarrows_\theta w$ iff one of the following holds:
 \begin{enumerate}
 \item $u=w$, $v\in \Sigma^+$ or $v\in P_\theta$.
     \item $u =\theta(w)$ and either $v\in \Sigma^+$ or $v\in P_\theta$ or $v=\gamma w$ for some $\gamma \in P_{\theta}$.
     \item $u=s^m =w$ and $v =s^n$, for $m,n\ge 1$ and either $s \in P_{\theta}$ or $s\in \Sigma^+$.
     \item $u =\theta(w)=(pq)^{j+1}p$, $v= (pq)^jp$ for some $p,q\in \Sigma^*$ and $j\ge 0$.
     \item $u =w=(pq)^{j+1}p$, $v= (pq)^jp$ for some $p,q\in P_{\theta}$ and $j\ge 0$.
     \item $u=w = xy$ and $v=(xy)^ix$, for $x,y \in P_{\theta}$ and $i\ge 0$.
      \item  $u = (xy)^{j+1}x=w$ and $v= yx$ for $x,y\in P_{\theta}$ and $j\ge 0$.
     \item $u=xy=\theta(w)$, $v= \theta(x)$ for $y\in P_{\theta}$.
     \item $u=xy=\theta(w)$, $v= x$ for $x, y\in P_{\theta}$.
     \item $u = \theta(t)s^{k}=\theta(w)$ and $v=s^{n}t $  where $s=t\theta(t)$.
 \end{enumerate}
\end{theorem}

\section{Solutions to $u\leftrightarrows_{\theta} L = L \leftrightarrows_{\theta} v$}\label{equations}

In this section we discuss solutions to the equation $u\leftrightarrows_{\theta} L = L \leftrightarrows_{\theta} v$ where $u,v \in \Sigma^+$ and $L\subseteq \Sigma^+$ which is a generalization of the equation $u\leftrightarrows_{\theta} w = w \leftrightarrows_{\theta} v$ where now $w$ is replaced with a set. Section \ref{concom} gave a complete characterization of words $u$ and $v$ when $L$ is a singleton. In this section we give solutions to the eqution $u\leftrightarrows_{\theta} L = L \leftrightarrows_{\theta} v$ under some special cases. 

We first recall the following from \cite{shyr1983} which characterizes languages such that $uL=Lv$ for non empty words $u$ and $v$. 

 \begin{proposition}\cite{shyr1983}\label{eqth1}
 Let $u,v\in \Sigma^+$ and $L\subseteq \Sigma^+$. Then $uL =Lv$ iff there exists $x,y \in \Sigma^*$ with  $|xy|\ge 1$ such that $u=(xy)^i$ and $v =(yx)^i$ for some  $i\ge 1$ and $L= \{x(yx)^j~:~ j\ge 0\}$.
\end{proposition}

  The following result gives solution to some simultaneous involution conjugate equations (\cite{Kari2010}).

  \begin{proposition}\cite{Kari2010}\label{eqth3}
      Let $x, y \in \Sigma^+$ and $\theta$ be an antimorphic involution with $xy = \theta(y)\theta(x)$ and $x\theta(y) = y\theta(x)$. Then,  $x = (\alpha\beta)^m$, $y = \alpha(\beta\alpha)^n$ with both $\alpha, \beta \in P_{\theta}$ for some $m \ge 1$ and $n \ge 0$.
  \end{proposition}

   We also recall the following results from \cite{Fan2011} which deals with some language equations incorporating the involution function.
 
\begin{proposition}\cite{Fan2011}\label{eqth2}
Let $\theta$ be an antimorphic involution, $u,v\in \Sigma^+$ and $L\subseteq \Sigma^+$.
If $\theta(L)u = vL$, then for $x\in \Sigma^+$, $y,z\in \Sigma^*$ with $xy\in \mathcal{Q}$, $v =(xy)^iz$, $u = z(\theta(y)\theta(x))^i$ for some $i\ge 1$ and  $$L \subseteq \{wz(\theta(y)\theta(x))^i~:~ w,z \in P_{\theta}, ~w\in \Sigma^*\}$$
\end{proposition}
We use the following lemma.
\begin{lemma}\label{reflem8}
For an antimorphic involution $\theta$, if either $xxy = yx\theta(x)$  or $ xxy= y\theta(x)x$  then,  $x$ and $y$ are powers of a common $\theta$-palindromic word. 
\end{lemma}
\begin{proof}
    We only prove the case when $xxy =yx\theta(x)$ as the proof for $xxy = y\theta(x)x$ is similar and we omit it. Let $xxy = yx\theta(x)$. Then by Lemma \ref{lemconj}, we have $xx =pq$, $y = (pq)^jp$ for some $j\ge 0$ and $x\theta(x) = qp$ for some $p,q\in \Sigma^+$. If $|x| \le |p|$ then, $x =p_1 = p_2q$ for $p=p_1p_2$. Then, $x\theta(x) = p_1\theta(q)\theta(p_2)= qp $ which implies that $p_2 \in P_\theta$ and $p_1\theta(q) = qp_1$ and by Lemma \ref{lemconj}, there exists $\alpha, \beta \in P_\theta$ such that $q= \alpha\beta$ and $p_1 = (\alpha\beta)^k\alpha$. Thus, we have $x = p_1 = (\alpha\beta)^k\alpha = p_2q = p_2\alpha\beta$ which implies that $\alpha\beta = \beta\alpha$ and by Lemma \ref{reflem7}, $\alpha$ and $\beta$ are powers of a common word. Hence, $x$ and $y$ are powers of a common $\theta$-palindromic word. The proof for the case when $|x|\ge |p|$ is similar.
    
\end{proof}
\begin{corollary}\label{refcor8}
For an antimorphic involution $\theta$, if either $x(xy)^i = (yx)^i\theta(x)$  or $ x(xy)^i = (y\theta(x))^ix$ for $i\ge 1$ then $x,y$ are powers of a common $\theta$-palindromic word. 
\end{corollary}
\begin{proof}
 We only prove for one of the given equation as the proof of the other one is similar.
    Given that $x(xy)^i = (yx)^i\theta(x)$. The case when $i=1$ is proved in Lemma \ref{reflem8}. Let $i\ge 2$. If $|x| \le |y|$ then $y=xp = q\theta(x)$ and by Proposition \ref{gg1} either $y= xs\theta(x)$ where $s \in \Sigma^*$ or $y= us\theta(u)$ where $x =us$ and $s\in P_\theta$. In both cases, $x(xy)^i = (yx)^i\theta(x)$ implies that, $xxs = s\theta(x)x$ and by Lemma \ref{reflem8}, both $s$, $x$ and hence, $y$ are powers of a common $\theta$-palindromic word. The case when $|x|\ge |y$ is similar and we omit it. 
\end{proof}
\begin{theorem}\label{teq1}
    Let $u\in \Sigma^+$ and $L\subseteq \Sigma^+$ such that $u\leftrightarrows_\theta L = L\leftrightarrows_\theta v$ and $uL =Lv$. Then, one of the following hold true.
    \begin{enumerate}
        \item $u = s^m$, $v = s^n$ and $L = \{s^k~:~k\ge 0\}$ for some $s\in \Sigma^+$.
        \item $u = s^m$, $v = s^n$ and $L = \{s^k~:~k\ge 0\}$ for some $s\in P_\theta$.
        \item $u=(xy)^i$, $v =(yx)^i$ for some  $i\ge 1$ and $L= \{x(yx)^j~:~ j\ge 0\}$ where $x,y\in P_\theta$.
        
    \end{enumerate}
\end{theorem}
\begin{proof}
    Given that $uL = Lv$ and by Proposition \ref{eqth1}, there exists $x,y \in \Sigma^*$ with  $|xy|\ge 1$ such that $u=(xy)^i$ and $v =(yx)^i$ for some  $i\ge 1$ and $L= \{x(yx)^j~:~ j\ge 0\}$.
    Since $u\leftrightarrows_\theta L = L\leftrightarrows_\theta u$, we have, 
    \begin{align*}
        u\leftrightarrows_\theta L &= 
        \{(xy)^{i+j}x, x(yx)^j(xy)^i,  (xy)^i\theta(x)(\theta(y)\theta(x))^j, \\
        & ~~~~~~\theta(x)(\theta(y)\theta(x))^j(xy)^i,  (\theta(y)\theta(x))^ix(yx)^j, \\
        &~~~~~x(yx)^j(\theta(y)\theta(x))^i, (\theta(y)\theta(x))^i\theta(x)(\theta(y)\theta(x))^{j},\\ &~~~~~~\theta(x)(\theta(y)\theta(x))^{j+i} ~:~ i\ge 1,j\ge 0\}
        \end{align*}
   \begin{align*}
        L\leftrightarrows_\theta v &= \{(xy)^{j+i}x, (yx)^i(xy)^jx,   (yx)^i\theta(x)(\theta(y)\theta(x))^j, \\
        & ~~~~~~\theta(x)(\theta(y)\theta(x))^j(yx)^i,  (\theta(x)\theta(y))^ix(yx)^j, \\
        &~~~~~x(yx)^j(\theta(x)\theta(y))^i, (\theta(x)\theta(y))^i\theta(x)(\theta(y)\theta(x))^{j},\\ &~~~~~~\theta(x)(\theta(y)\theta(x))^{j}(\theta(x)\theta(y))^i ~:~ i\ge 1,j\ge 0\}
        \end{align*}
    We now have the following cases.
    \begin{enumerate}
        \item Let $(xy)^{j}x(xy)^i = (yx)^mx(yx)^n$ where $i+j =m+n$, $i,m\ge 1$.  If $j \neq 0$ then, $xy =yx$ which implies by Lemma \ref{lpow} that both $x$ and $y$ are powers of a common word. 
        If  $j =0$  then, $xxy = yxx$ and by Lemma \ref{reflem7}, $x$ and $y$ are powers of a common word. Hence in both cases, $u = s^m$, $v = s^n$ and $L = \{s^k~:~k\ge 0\}$ for some $s\in \Sigma^+$.
        \item  Let $(xy)^{j}x(xy)^i = (yx)^m\theta(x)(\theta(y)\theta(x))^n$ where $i+j =m+n$, $i,m\ge 1$. We now have the following subcases. 
        \begin{itemize}
            \item If $j\neq  0$ then, $xy=yx$ which implies by Lemma \ref{lpow} that both $x$ and $y$ are powers of a common word $s\in \Sigma^+$. If in addition $n\neq 0$ then, $xy \in P_{\theta}$ which implies that $s\in P_\theta$. If $n=0$ then, $x\in P_\theta$ and hence, both $x$ and $y$ are powers of a common word $s\in \Sigma^+$.
            \item If both $j=0$ and $n=0$ then, $i=m$ and $x(xy)^i= (yx)^i\theta(x)$
            by Corollary \ref{refcor8}, $x$ and $y$ are powers of a common $\theta$-palindromic word.
            \item If $j= 0$ and $n\neq 0$ then, both $x(xy)^i = (yx)^m\theta(x)(\theta(y)\theta(x))^n$ which implies that $xy \in P_\theta$ and $x(xy)^i =(yx)^i\theta(x)$. Then, by Corollary \ref{refcor8}, $x$ and $y$ are powers of a common $\theta$-palindromic word.           \end{itemize}
            \item Let $(xy)^{j}x(xy)^i = \theta(x)(\theta(y)\theta(x))^n(yx)^m$ where $i+j =m+n$, $i,m\ge 1$. Then, $xy=yx$ and $x\in P_\theta$ which implies by Lemma \ref{lpow} that both $x$ and $y$ are powers of a common $\theta$-palindromic word. 
            \item Let $(xy)^{j}x(xy)^i = (\theta(x)\theta(y))^mx(yx)^n,$ where $i+j =m+n$, $i,m\ge 1$. We now have the following subcases. 
        \begin{itemize}
            \item If $j\neq 0$ then, $x,y\in P_\theta$.  If in addition $n=0$ then $x(xy)^i = (\theta(x)\theta(y))^ix = (xy)^ix$ and if $n\neq 0$, we also get $xy=yx$. Hence, by Lemma \ref{lpow} both $x$ and $y$ are powers of a common $\theta$-palindromic word.
            \item If $j=0$ then, $x\in P_\theta$ and $xy = yx$.  Hence, by Lemma \ref{lpow} both $x$ and $y$ are powers of a common $\theta$-palindromic word.
        \end{itemize}
         \item Let $(xy)^{j}x(xy)^i = x(yx)^n(\theta(x)\theta(y))^m,$ where $i+j =m+n$, $i,m\ge 1$ which implies that both $x,y\in P_\theta$. If $i\neq m$ then, $xy =yx$ and by Lemma \ref{lpow} both $x$ and $y$ are powers of a common $\theta$-palindromic word.
         \item Let $(xy)^{j}x(xy)^i = \theta(x)(\theta(y)\theta(x))^{m+n}$, where $i+j =m+n$, $i,m\ge 1$. Then, $xy \in P_\theta$ and we have the following subcases. 
        \begin{itemize}
            \item If $j=0$ then, $x(xy)^i = \theta(x)(\theta(y)\theta(x))^{m+n}$ which implies that $x \in P_\theta$ and $xy = \theta(y)x$. Hence, by Lemma \ref{lemconj}, $y = qp$, $x=(pq)^tp$ for $p,q\in P_\theta$. Then,
            \begin{align*}
        u\leftrightarrows_\theta L &= 
        \{[(pq)^{t+1}p]^{i+j}(pq)^tp, [(pq)^{t+1}p]^{j}(pq)^tp[(pq)^{t+1}p]^{i}, \\
        & ~~~~~~[(pq)^{t+1}p]^{i}(pq)^tp[(pq)^{t+1}p]^{j},  (pq)^tp[(pq)^{t+1}p]^{i+j}, \\
        &~~~~~~:~ i\ge 1,j\ge 0\}
        \end{align*}
        Also, observe that $[(qp)(pq)^tp]^i[(pq)^{t+1}p]^j(pq)^tp \in L\leftrightarrows_\theta v$.
        Since, $u\leftrightarrows_\theta L = L\leftrightarrows_\theta v$, we have either $qp=pq$ or $qpp =ppq$, which implies that both $p$ and $q$ are powers of a common word. 
            \item If $j\neq 0$ and $n\neq 0$ then, $x, y, xy \in P_\theta$ and by Lemma \ref{lpow}, both $x$ and $y$ are powers of a common $\theta$-palindromic word.
        \end{itemize}
        \item Let $(xy)^{j}x(xy)^i = \theta(x)(\theta(y)\theta(x))^{n}(\theta(x)\theta(y))^{m}$, where $i+j =m+n$, $i,m\ge 1$. Then, $x,y\in P_\theta$. If $i\neq m$  then, $xy=yx$ and by Lemma \ref{lpow}, both $x$ and $y$ are powers of a common $\theta$-palindromic word.
    \end{enumerate}
    Hence, the result. 
    \end{proof}
 We now have the following which follows directly from Proposition \ref{eqth1} and Theorem \ref{teq1}.   
\begin{theorem}\label{teq2}
    Let $u\in \Sigma^+$ and $L\subseteq \Sigma^+$ such that $u\leftrightarrows_\theta L = L\leftrightarrows_\theta v$ and $uL =L\theta(v)$. Then, one of the following hold true.
    \begin{enumerate}
        \item $u = s^m$, $v = \theta(s)^n$ and $L = \{s^k~:~k\ge 0\}$ for some $s\in \Sigma^+$.
        \item $u = s^m$, $v = s^n$ and $L = \{s^k~:~k\ge 0\}$ for some $s\in P_\theta$.
        \item $u=(xy)^i= v$ for some  $i\ge 1$ and $L= \{x(yx)^j~:~ j\ge 0\}$ where $x,y\in P_\theta$.
        \end{enumerate}
\end{theorem}

One can easily observe from Remark \ref{rem1}, the following.
$$u\leftrightarrows_\theta L = L_1\leftrightarrows_\theta u_1=u_1\leftrightarrows_\theta L_1$$
for $u_1\in u_\theta$ and $L_1\in L_\theta$. Hence by Proposition \ref{eqth2} we conclude the following.

\begin{theorem}\label{teq3}
    Let $u\in \Sigma^+$ and $L\subseteq \Sigma^+$ such that $u\leftrightarrows_\theta L = L\leftrightarrows_\theta v$ and $uL =\theta(L)v$. Then, for $x\in \Sigma^+$, $y,z\in \Sigma^*$ with $xy\in \mathcal{Q}$, $u =(xy)^iz$, $v = z(\theta(y)\theta(x))^i$ for some $i\ge 1$ and  $$L \subseteq \{wz(\theta(y)\theta(x))^i~:~ w,z \in P_{\theta}, ~w\in \Sigma^*, ~i\ge 1\}$$
\end{theorem}
\begin{proof}
Observe that,
$$ u\leftrightarrows_\theta L = \{uL, Lu, L\theta(u), \theta(u)L, u\theta(L), \theta(L)u, \theta(u)\theta(L), \theta(L)\theta(u)\}$$
$$ L\leftrightarrows_\theta v = \{vL, Lv, L\theta(v), \theta(v)L, v\theta(L), \theta(L)v, \theta(v)\theta(L), \theta(L)\theta(v)\}$$
    Given that $uL = \theta(L)v$, which implies that $\theta(v)L = \theta(L)\theta(u)$ and by Proposition \ref{eqth2}, $u =\theta(v)$. Hence, 
   $$\{ Lu, L\theta(u), \theta(u)L, u\theta(L), \theta(L)u, \theta(u)\theta(L), \}$$
 $$  = \{\theta(u)L, L\theta(u), Lu,  \theta(u)\theta(L), u\theta(L), \theta(L)u\}$$
   Thus, by Proposition \ref{eqth2}, for $x\in \Sigma^+$, $y,z\in \Sigma^*$ with $xy\in \mathcal{Q}$, $u =(xy)^iz$, $v = z(\theta(y)\theta(x))^i$ for some $i\ge 1$ and  $$L \subseteq \{wz(\theta(y)\theta(x))^i~:~ w,z \in P_{\theta}, ~w\in \Sigma^*, ~i\ge 1\}$$
    
\end{proof}

 We now have the following which follows directly from Proposition \ref{eqth2} and Theorem \ref{teq3}. 

\begin{theorem}\label{teq4}
    Let $u\in \Sigma^+$ and $L\subseteq \Sigma^+$ such that $u\leftrightarrows_\theta L = L\leftrightarrows_\theta v$ and $uL =\theta(L)\theta(v)$. Then, for $x\in \Sigma^+$, $y,z\in \Sigma^*$ with $xy\in \mathcal{Q}$, $u =(xy)^iz= v$ for some $i\ge 1$ and  $$L \subseteq \{wz(\theta(y)\theta(x))^i~:~ w,z \in P_{\theta}, ~w\in \Sigma^*, ~i\ge 1\}$$
\end{theorem}

\begin{lemma}\label{lemeq5}
    Let $u,v\in \Sigma^+$ and $L\subseteq \Sigma^+$. The following are true. 
    \begin{enumerate}
        \item If $uL =vL$ then $u=v$. 
        \item If $uL = \theta(v)L$ then $u=\theta(v)$.
        \item If $uL = v\theta(L)$. then $u=v$ and $L=\theta(L)$.
         \item If $uL = \theta(v)\theta(L)$ then, $u=\theta(v)$ and $L=\theta(L)$.
    \end{enumerate}
 \end{lemma}
 \begin{proof}
     We only prove the first implication as the others are similar.
     Let $w\in L$ be such that $|w|\le |x|$ for all $x\in L$. Since $uL = vL$, $uw = vx$
 for some $x\in L$ and $|x|\ge |w|$ which implies that $|u|\ge |v|$. 
 Also, there exists a $y\in L$ such that $|w|\le |y|$ and $uy = vw$.
 If $|x|>|w|$ then, $|u|>|v|$ and hence, $|wv| <|uy|$ a contradiction. Similarly, we can show that $|x|\nless |w|$. Hence, $|w| =|x|$ which implies that $|u|=|v|$. Since, $uL=vL$ we conclude that $u=v$.  \end{proof}
 By Lemma \ref{lemeq5}, we conclude the following.
 \begin{theorem}\label{teq5}
      Let $u\in \Sigma^+$ and $L\subseteq \Sigma^+$ such that $u\leftrightarrows_\theta L = L\leftrightarrows_\theta v$. The following are true.
      \begin{enumerate}
        \item If $uL =vL$ then $u=v$. 
        \item If $uL = \theta(v)L$ then $u=\theta(v)$.
        \item If $uL = v\theta(L)$. then $u=v$ and $L=\theta(L)$.
         \item If $uL = \theta(v)\theta(L)$ then, $u=\theta(v)$ and $L=\theta(L)$.
    \end{enumerate}
 \end{theorem}

 \section{Conclusions}
\label{conclusions}
This paper defines and investigates the binary word operation strong-$\phi$-bi-catenation which, when iteratively applied to words $u$ and $v$ generates words in the set $\{u,\phi(u),v,\phi(v)\}^+$.
The operation was naturally extended to languages (Section \ref{Lext}) and we investigated some of its properties.
 Future topics of research include extending the  $\leftrightarrows_\phi$-conjugacy on words (Section \ref{concom} ) to $\leftrightarrows_\phi$-conjugacy on languages as well as exploring an associative version of strong-$\theta$-bi-catenation.

\bibliographystyle{te}

\begin{thebibliography}{1}

\bibitem{Czeizler10}Czeizler, E., Kari, L. \& Seki, S. On a special class of primitive words. {\em Theoretical Computer Science}. \textbf{411}, 617-630 (2010)
\bibitem{shyr1983}Shyr, H. On two languages that commute. {\em Notes On Semigroups}. \textbf{9} pp. 257-269 (1983)
\bibitem{Bcat}Shyr, H. \& Yu, S. Bi-catenation and shuffle product of languages. {\em Acta Informatica}. \textbf{35} pp. 689-707 (1998)
\bibitem{Fan2011}Fan, C. \& Huang, C. Solutions to the involution language equation. {\em International Journal Of Computer Mathematics}. \textbf{88-11} pp. 2285-2292 (2011)
\bibitem{Seki10}Kari, L. \& Seki, S. An improved bound for an extension of Fine and Wilf’s theorem and its optimality. {\em Fundamenta Informaticae}. \textbf{101} pp. 215-236 (2010)
\bibitem{twin}Kari, L., Mahalingam, K. \& Seki, S. Twin Roots and their Properties. {\em Theoretical Computer Science}. \textbf{410}, 2393-2400 (2009)
\bibitem{masson}Kari, L., Masson, B. \& Seki, S. Properties of Pseudo-Primitive Words and Their Applications. {\em International Journal Of Foundations Of Computer Science}. \textbf{22}, 447-471 (2011)
\bibitem{Knuth77}Knuth, D., Morris, J. \& Pratt, V. Fast pattern matching in strings. {\em SIAM Journal On Computing}. \textbf{6}, 323-350 (1977)
\bibitem{deLuca81}Luca, A. A combinatorial property of the Fibonacci words. {\em Information Processing Letters}. \textbf{12}, 193-195 (1981)
\bibitem{Seebold85}Séébold, P. Sequences generated by infinitely iterated morphisms. {\em Discrete Applied Mathematics}. \textbf{11}, 255-264 (1985)
\bibitem{Shallit88}Shallit, J. A generalization of automatic sequences. {\em Theoretical Computer Science}. \textbf{61}, 1-16 (1988)
\bibitem{Lothaire97}Lothaire, M. Combinatorics on words. (Cambridge University Press,1997)
\bibitem{hsiao}Hsiao, H., Huang, C. \& Yu, S. Word operation closure and primitivity of languages. {\em Journal Of Automata Languages And Combinatorics}. \textbf{19}, 157-171 (2014)
\bibitem{kinvcat}L.Kari \& K.Mahalingam k-involution codes and related sets. {\em Journal Of Discrete Mathematical Sciences And Cryptography}. \textbf{10} pp. 485-503 (2007)
\bibitem{man14}Kari, L. \& Kulkarni, M. Generating the Pseudo-powers of a word. {\em Journal Of Universal Computer Science}. \textbf{8}, 243-256 (2014)
\bibitem{Yu00}Yu, S. \& Zhao, Y. Properties of Fibonacci languages. {\em Discrete Mathematics}. \textbf{224}, 215-223 (2000)
\bibitem{Kari2010}Kari, L. \& Mahalingam, K. Watson-Crick Palindromes in DNA Computing. {\em Natural Computing}. \textbf{9}, 297-316 (2010)
\bibitem{Kari08}Kari, L. \& Mahalingam, K. Watson-Crick Conjugate and Commutative Words. {\em DNA Computing}. \textbf{4848} pp. 273-283 (2008)
\bibitem{jonoska2004}Jonoska, N. \& Mahalingam, K. Languages of DNA based code words. {\em DNA Computing}. \textbf{2943} pp. 61-73 (2004)
\bibitem{tulpan2003}Tulpan, D., Hoos, H. \& Condon, A. Stochastic local search algorithms for DNA word design. {\em DNA Computing}. \textbf{2568} pp. 229-241 (2003)
\bibitem{mauri2004}Mauri, G. \& Ferretti, C. Word design for molecular computing: A survey. {\em DNA Computing}. \textbf{2943} pp. 37-47 (2004)
\bibitem{cai2023}Kari, L. \& K.Mahalingam Watson-Crick Powers of a word. {\em  Algebraic Informatics}. \textbf{13706} pp. 136-148 (2022)
\bibitem{Lshuf}M.Ito, L. \& G.thierrin Shuffle and scattered deletion closure of languages. {\em Theoretical Computer Science}. \textbf{245-1} pp. 115-133 (2000)
\bibitem{kcat}L.Kari \& G.thierrin K-catenation and applications: k-prefix codes. {\em Journal Of Information And Optimization Sciences}. \textbf{16-2} pp. 263-276 (1995)
\bibitem{lyndon}Lyndon, R. \& Schützenberger, M. The equation a<sup>M</sup>=b<sup>N</sup>c<sup>P</sup> in a free group. {\em Michigan Mathematical Journal}. \textbf{9} pp. 289-298 (1962)
\bibitem{Wu92}H. L. Wu On the Properties of Primitive Words. (Institute of Applied Mathematics, Chung-Yuan Christian University,1992)
\bibitem{Kari02}Kari, L., Kitto, R. \& Thierrin, G. Codes, Involutions, and DNA Encodings. {\em Formal And Natural Computing: Essays Dedicated To Grzegorz Rozenberg}. \textbf{2300} pp. 376-393 (2002)
\bibitem{Lthesis}L.Kari On Insertion and Deletion in Formal Languages. (University of Turku,1991)
\bibitem{ssyu2005}Yu, S. Languages and Codes. (Tsang Hai Book Publishing Co.,2005)
\bibitem{Eilen}Eilenberg, B. Automata, Languages and Machines. (Academic Press,1974)
\bibitem{paun1998}Păun, G., Rozenberg, G. \& Salomaa, A. DNA Computing: New Computing Paradigms. {\em Texts In Theoretical Computer Science}. (1998)
\end{thebibliography}

\vspace{1cc}

\noindent
{\bf Kalpana Mahalingam \\
Department Of Mathematics, \\
Indian Institute of Technology, \\
Chennai-600036. India.\\
Email~:~{kmahalingam@iitm.ac.in}
\end{document}